\documentclass[11pt, letterpaper]{article}
\usepackage[affil-it]{authblk}
\usepackage[in]{fullpage}
\usepackage{graphicx}
\usepackage{graphicx}
\usepackage{amsfonts}
\usepackage{amsmath, amsthm, amssymb}
\usepackage[numbers,sort&compress]{natbib}
\usepackage{url}
\usepackage{color}
\usepackage{multirow}
\usepackage{enumerate}
\usepackage{hyperref}
\usepackage{epstopdf}
\usepackage{tikz}
\usepackage{here}
\usepackage{bigstrut}
\usepackage{doi}
\newcommand{\trace}{\textbf{trace}}
\newcommand{\diag}{\textbf{diag}}
\newcommand{\Diag}{\textbf{Diag}}
\newcommand{\bdw}{\textbf{bdw}}
\newcommand{\vect}{\textbf{vec}}

\newcommand{\minpart}{\operatornamewithlimits{minPart}}
\newcommand{\proofstep}[1]{\textbf{Step #1:}}

\newtheorem{thm}{Theorem}
\newtheorem{lem}[thm]{Lemma}
\newtheorem{cor}[thm]{Corollary}
\newtheorem{prop}[thm]{Proposition}

\newtheorem{example}{Example}

\newcommand{\R}{\mathbb{R}}

\newcommand{\cS}{\mathcal{S}}
\newcommand{\transpose}[1]{\ensuremath{#1^{\mathsf{T}}}}

\usepackage{xspace}
\xspaceaddexceptions{[]\{\}}
\newcommand{\BandwidthProblem}{Bandwidth Problem\xspace}

\def \N{\mathbb{N}}

\renewcommand{\binom}[2]{\genfrac(){0pt}{}{#1}{#2}}

\providecommand{\keywords}[1]{\textbf{\textit{Keywords.}} #1}

\title{Lower Bounds for the Bandwidth Problem
\footnote{ This project has received funding from the European Union's
Horizon 2020 research and innovation programme under the Marie Sk\l{}odowska-Curie
grant agreement No 764759 and the Austrian Science Fund (FWF): P 28008-N35.}
}
\author[1]{Franz Rendl}
\author[2]{Renata Sotirov}
\author[1]{Christian Truden}

\affil[1]{Department of Mathematics, Alpen-Adria-Universit\"at, Klagenfurt, Austria}
\affil[2]{Department of Econometrics and OR, Tilburg University, The Netherlands}

\date{}

\begin{document}

\maketitle

\begin{abstract}
  The \BandwidthProblem  seeks for a   simultaneous permutation of the rows and
  columns of the  adjacency matrix of a graph such that all nonzero entries are
  as close as possible to the main diagonal.
  This work focuses on investigating novel approaches  to obtain lower bounds
  for the bandwidth problem.
  In particular, we use   vertex partitions  to  bound the bandwidth of a graph.
  Our approach contains prior approaches for bounding the bandwidth  as special cases.
  By varying sizes of partitions, we achieve a trade-off between quality of bounds
  and efficiency of computing them.
  To compute lower bounds, we  derive a Semidefinite Programming relaxation.
  We evaluate the performance of our approach on several data sets, including real-world instances.
  \\[1em]
  \keywords{
  Bandwidth Problem, Graph Partition, Semidefinite Programming.}
\end{abstract}

\section{Introduction}

The Bandwidth Problem (BP) is the problem of labeling the vertices of a given
undirected graph with distinct integers such that the maximum
difference between
the labels of adjacent vertices is minimal.
It originated in the 1950s from sparse matrix computations,
and received much attention since
Harary's \cite{Harary:67} description of the problem and Harper's paper
\cite{Harper66} on the bandwidth of the $n$-cube (see also \cite{ChChDeGi:82,DiPeSe:02}).
Berger-Wolf and Reingold \cite{BergWoRe:02} showed that the
problem of designing a code to minimize distortion in multi-channel
transmission can be formulated as the Bandwidth Problem for  generalized
Hamming graphs.
The BP belongs to a class of combinatorial optimization problems known as graph layout problems.
The Cyclic Bandwidth \cite{RODRIGUEZTELLO201517,dKlSoNa}, Cutwidth \cite{MARTI2013137,CAVERO2021105116},
Antibandwidth \cite{1fb351e298a049cc855e645215159c30} and Linear Arrangement Problem \cite{Harper1964,RODRIGUEZTELLO20083331} also belong to this class of problems.
The Bandwidth Problem arises in many different engineering
applications related
to efficient storage and processing. It also plays a role in designing
parallel computation networks, VLSI layouts, and constraint
satisfaction problems,
see e.g.,~\cite{ChChDeGi:82,DiPeSe:02,LaWi:99} and the references therein.

Determining the bandwidth is
NP-hard \cite{Papadem:76} and even approximating the bandwidth within a given
factor is known to be  NP-hard \cite{Unger1998}.
Moreover, the BP is known to be NP-hard even on trees with maximum degree three
\cite{GagrJoKn:78} and on caterpillars with hair length three \cite{Monien}.
On the other hand, the Bandwidth Problem has been solved for a few families of
graphs having special properties. Among these are the path,
the complete graph, the complete bipartite graph \cite{Chvatal1},
the hypercube
graph \cite{Harper66}, the grid graph
\cite{Chvalatalova}, the complete $k$-level $t$-ary tree \cite{Smith}, the
triangular graph \cite{HwLa77}, and the
triangulated triangle \cite{HoMcDSa}.
Blum et al.~\cite{BlKoRaVe} and
Dunagan and Vempala \cite{DuVe:01}
propose an $O(\log^3 n \sqrt{\log \log n })$ approximation
algorithm for the
bandwidth, where $n$ is the number of vertices.

Several lower and upper bounding approaches for the bandwidth of a
graph are considered in the literature.
Cuthill and McKee \cite{CutMcKee} proposed a heuristic to relabel
the vertices of the graph so as to reduce the bandwidth after
relabeling. It is widely used in practice, see for instance
\cite{Turner86}.
MATLAB offers
the command \texttt{symrcm}
as an implementation of this heuristic.
For graphs with symmetry there exists  an improved reverse Cuthill-McKee
algorithm, see \cite{DamSotirov15}.
However, it is much more difficult to obtain lower bounds on the bandwidth.
The following two approaches have been proposed in the literature.

\paragraph{Lower bounds based on 3-partitions}
Juvan and Mohar \cite{JuMo:99} consider 3-partitions of the vertices
into partition blocks $S_{1}, S_{2}, S_{3}$
of (fixed) sizes $m_{1},m_{2}$ and $m_{3}$.
If all such partitions have edges joining $S_{1}$ and $S_{3}$,
then clearly the bandwidth must be bigger than $m_{2}$.
Juvan and Mohar introduce eigenvalue-based lower bounds on the bandwidth
which were refined by
Helmberg et al.~\cite{HelReMoPo:95}  leading to the following
bound based on eigenvalues of the Laplacian $L$ of the graph
$$
{\bdw} >\frac{ n \lambda_{2}(L)}{\lambda_{n}(L)},
$$
see also the subsequent section.
The same lower bound was derived by
Haemers \cite{haemers} by exploiting interlacing of
Laplacian eigenvalues.
Povh and Rendl  \cite{PoRe:07} showed that this eigenvalue bound
can also be obtained by solving
a Semidefinite Programming (SDP) relaxation for
a special Minimum Cut (MC) problem.
They further tightened the SDP relaxation and consequently obtained a stronger
lower bound for the Bandwidth Problem.
Rendl and Sotirov \cite{Rendl2017} showed
how to further tighten the SDP
relaxation from \cite{PoRe:07}.

\paragraph{Bounds based on permutations}
A labeling of the vertices of a graph corresponds to a simultaneous
permutation of the rows and columns of the adjacency matrix.
This may be expressed by pre- and post-multiplication with a
permutation matrix, leading to quadratic assignment formulations
of the bandwidth.
De Klerk et al.~\cite{dKlSoNa} proposed two lower
bounds
based on SDP relaxations of
the resulting Quadratic Assignment Problem (QAP).
The numerical results in \cite{dKlSoNa} show that both
their bounds dominate the
bound of Blum et al.~\cite{BlKoRaVe}, and that in most of
the cases their bounds  are stronger than the bound by
Povh and Rendl \cite{PoRe:07}.

In \cite{DamSotirov15}, the authors derived an SDP relaxation
of the minimum cut
problem by strengthening the well-known SDP relaxation for the QAP.
They derive strong bounds for  the bandwidth  of highly symmetric
graphs with up
to 216 vertices by exploiting symmetry.
For general graphs, their approach is rather restricted.
Above mentioned bounds are either unsatisfyingly  weak, or
computing them is challenging already for small (general) graphs, i.e., graphs of
about $30$ vertices.

\paragraph{Our contribution}
We introduce a general $k$-partition model to get lower bounds on the
bandwidth. It contains (with $k=3$) the 3-partition model from Juvan and
Mohar \cite{JuMo:99} and (with $k=n$) the permutation-based formulation of
the problem, see Section \ref{sect:modelApproach} below.
The $k$-partition problem is still NP-complete.
Therefore, we introduce tractable relaxations
based on SDP.
In Section \ref{sect:SDPmodels} such a relaxation based
on the ``matrix-lifting'' idea  is introduced. It leads to
an SDP in matrices of order $n \cdot k$.
It is known that the feasible region of
such a relaxation always has a nullspace of dimension $n+k-1$.
We identify an $n$-dimensional part of this nullspace, which can be eliminated
using a simple combinatorial argument.
Finally, in Section \ref{sect:numerResults}, we show
that the new partition model leads to improved lower
bounds for the bandwidth, even in case of small values of $k$, like $k\leq 6$.
Moreover, we  provide strong bounds for graphs with up to $128$ vertices
in a reasonable time frame.

\paragraph{Notation}
The space of $n\times n$ symmetric matrices is denoted by $\cS_n$ and the space of $n\times n$ symmetric positive
semidefinite  matrices by $\cS^+_n$. For two matrices $X,Y\in \R^{n\times n}$, $X\geq Y$, means $x_{ij}\geq y_{ij}$, for all $i,j$.
The set of $n \times n$ permutation matrices is denoted by $\Pi_n$.
Further, for a matrix $A$ the corresponding transposed matrix is denoted by $\transpose{A}$ while
$A^{\bot}$ denotes the orthogonal complement.
We use $I_n$ to denote the identity matrix of order $n$, and $e_n^i$ to denote the  $i$-th standard basis vector of length $n$.
Similarly, $J_n$ and $e _n$ denote the $n \times n$ all-ones matrix and all-ones $n$-vector, respectively.

The trace operator is denoted by
$\trace$, and $\langle \cdot, \cdot\rangle$ denotes the trace inner product.
The Hadamard product of two matrices $A$ and $B$ of the same size
is denoted  by $A\circ B$ and defined as $(A\circ B)_{ij} = a_{ij}\cdot b_{ij}$ for all $i,j$.
The $\diag$ operator maps an $n\times n$ matrix to the $n$-vector given by its diagonal, while the $\vect$ operator
stacks the columns of a matrix in a vector. We denote by $\Diag$ the adjoint operator of $\diag$.

\section{The Bandwidth Problem} \label{sect:bwd}

We now formally introduce the Bandwidth Problem  as a
Quadratic Assignment Problem  with special data matrices $A$
and $B_{r,n}$.

Let $G=(V,E)$ be an undirected simple graph with $| V| =n$ vertices and edge set $E$.
A bijection $\phi : V  \rightarrow \left\{1,\ldots, n \right\} $
is called a \textit{labeling} of the vertices of $G$.
The bandwidth of a graph $G$ with respect to the labeling $\phi$ is defined as follows
$$
\bdw(\phi,G):=
\max_{[i,j] \in E} |  \phi(i)  - \phi (j) |.
$$
The \textit{bandwidth of a graph}  $G$ is defined as the minimum of
$\bdw(\phi, G)$ over all labelings $\phi$, i.e.,
\begin{align*}
  \bdw
  (G):= \min \left\{
  \bdw(\phi,G):   ~ \phi ~ \text{labeling of~} G
  \right\}.
\end{align*}
Equivalently, one can   consider the adjacency matrix $A$ of the graph $G$.
The bandwidth of $A$  amounts to a simultaneous permutation of the rows and columns of the adjacency matrix  such that the largest distance of a
nonzero entry from the main diagonal is as small as possible. The \textit{bandwidth of an adjacency matrix} $A$ is defined as:
$$
\bdw(A):= \bdw(G).
$$
Therefore, from now on we assume that a graph  $G$ is given through its adjacency matrix $A$.
Since in terms of matrices the BP  seeks for a simultaneous permutation of the rows and columns of
$A$ such that all nonzero entries are as close as possible to the main diagonal, a ``natural'' problem formulation is as follows.

Let $r$ be an integer such that $1\leq r\leq n-2$, and $B_{r,n}=(b_{ij})$ be the symmetric matrix of order $n$ defined as follows
\[
b_{ij}:=
\begin{cases}
  1,  & \text{for } |i-j| > r, \\
  0, & \text{otherwise}.
\end{cases}
\]
Then, the following holds:
\begin{equation}
  \min\limits_{Q \in \Pi_n}  \langle \transpose{Q} A Q, B_{r,n} \rangle
  =
  \begin{cases}
    0,  		& \text{{if }} \bdw(A) \leq r,\\
    >	0,   & \text{{if }} \bdw(A) > r.
  \end{cases}
  \label{eq:qap}
\end{equation}
The  minimization problem has the form of a QAP, which might be even harder to solve than actually computing $\bdw(A)$.
The idea of formulating  the Bandwidth Problem as a QAP  was  suggested by Helmberg et al.~\cite{HelReMoPo:95}.
De Klerk et al.~\cite{dKlSoNa} considered two SDP-based bounds for the Bandwidth Problem that are obtained from the SDP relaxations for the QAP introduced in \cite{ZhKaReWo:98} and \cite{dKSo:12}.
The results show that it is hard to obtain bounds for  graphs with $32$ vertices, even though the symmetry in the graphs under consideration was exploited.

Since  it is very difficult to solve QAPs in practice for sizes larger than $30$ vertices other approaches are needed for deriving bounds for the bandwidth of graphs.

\section{Partition Approach} \label{sect:modelApproach}

We show how to use vertex partitions in order
to obtain lower bounds for the bandwidth of a graph.
For $3 \leq k \leq n$ let $m \in \N^{k}$ be given with $m_{i} \geq 1$ ($i=1,\ldots, k$), $\sum_{i=1}^{k} m_{i} =n$.
We consider partitions of the vertex set  $V$ into $k$ subsets $\{ S_{1}, \ldots, S_{k} \}$ such that $|S_{j}| = m_{j}, ~ j=1,\ldots,k$.
These are in one-to-one correspondence with $n \times k$ partition matrices:
\begin{equation} \label{PM}
  \mathcal{P}_{m} := \{X \in \{0,1\}^{n \times k}: Xe_{k}=e_{n},~ \transpose{X}e_{n}=m\},
\end{equation}
where for the partition $(S_{1}, \ldots, S_{k})$ we set $x_{ij}=1$
whenever $i \in S_{j}$, $i=1,\ldots, n$. Since any vertex $i \in V$ is assigned to precisely one
of the blocks $S_{j}$ we can define
the map $p:V=\{ 1,\ldots,n \} \mapsto \{1,\ldots,k\}$
given by
$$
p(i)=j \Leftrightarrow x_{ij}= 1 \Leftrightarrow  i \in S_{j},
$$
which identifies the partition block containing vertex $i$.
Thus, given the partition matrix $X \in \mathcal{P}_{m}$ we get
$S_{j} = \{i \in V: p(i)=j \}$
for all $1 \leq j \leq k$.
For $1 \leq r \leq k-2$ let $B_{r,k}=(b_{ij})$ be the 0--1 matrix of order $k$
with
\begin{equation} \label{Brk}
  b_{ij} = \begin{cases}
  1, & |i-j| > r, \\
  0, & |i-j| \leq r.
\end{cases}
\end{equation}
Suppose that $i \in S_{u}, j \in S_{v}$, i.e.,
$p(i)=u, ~p(j)=v$. Then for $X\in \mathcal{P}_{m}$ the following holds:
$$
(X B_{r,k} \transpose{X})_{ij} =
\transpose{{e_{k}^{u}}}  B_{r,k} e_{k}^{v}
=
\begin{cases}
  1, & |u-v| >r, \\
  0, & |u-v| \leq r.
\end{cases}
$$
Therefore we get
$$
\frac{1}{2}\langle A, X B_{r,k} \transpose{X} \rangle =
\sum_{\substack{i,j \in V,\\ i<j }  } a_{ij}(X B_{r,k} \transpose{X} )_{ij} =
\sum_{[i,j] \in E} (X B_{r,k} \transpose{X})_{ij} =
\sum_{\substack{[i,j] \in E, \\  |p(i) -p(j)|>r} } 1.
$$
Hence, this term counts the number of edges with endpoints
in partition blocks
of distance greater than $r$ under the map  $p$.

\paragraph{Basic Partition}
It will be convenient to introduce
the special partition matrix $\overline{X}$ corresponding to the
\textit{basic partition} $\overline{p}$ which assigns the first
$m_1$ vertices to $S_1$ the next $m_2$ vertices to $S_2$ and so on.
Thus, the $n \times k$ matrix
$\overline{X}$ is characterized by columns of consecutive blocks of ones
of appropriate lengths.
Therefore the $n \times n$ matrix
$$
B := \overline{X}B_{r,k}\transpose{\overline{X}}
$$
is a block matrix with blocks of sizes $m_i \times m_j$.
The nonzero blocks of this matrix correspond to all-ones matrices
of size $m_i \times m_j$ whenever the entry $(B_{r,k})_{ij}=1$, see also
Figure   \ref{fig_upperLowerBound} below.
Thus, for a given $n \times n$ adjacency matrix $A$ the term
$\frac{1}{2} \langle A,\overline{X} B_{r,k} \transpose{\overline{X}} \rangle$
counts
the number of
edges  joining vertices in partition blocks
of distance greater than $r$.

\paragraph{General Partition}
In general, any partition matrix $X  \in \mathcal{P}_{m}$ can be obtained
from the basic partition matrix  $\overline{X}$
by row-permutations that are defined by a permutation
matrix  $P \in \Pi_n$.
Thus
$$\mathcal{P}_m = \{ P \overline{X} \colon P \in \Pi_n  \},$$
where $\overline{X}$ is the basic partition matrix.
The following transformation is obtained by replacing $X$ by $P\overline{X}$:
\[
\frac{1}{2} 	\langle A, X B_{r,k} \transpose{X} \rangle =  \frac{1}{2}
\langle  A, P \overline{X} B_{r,k} \transpose{\overline{X}} \transpose{P} \rangle
=\frac{1}{2} 	\langle \transpose{P} A P,
\overline{X} B_{r,k} \transpose{\overline{X}}  \rangle.
\]
This shows that the permutation $P\in \Pi_n$ can be
applied either to  the adjacency matrix $A$ or to the
matrix $\overline{X}B_{r,k} \transpose{\overline X}$.

The following example serves as an illustration of this property.
\begin{example}
  We consider a $15\times 15$ matrix and the partitioning
  $m=\transpose{(3,3,3,\allowbreak 3,3)}$. Moreover, we choose	$r=2$.
  If
  $\langle  A,
  \overline{X} B_{r,k} \transpose{\overline{X}}  \rangle >0$,
  then there must be an
  edge with endpoints	in blocks of distance larger than $r=2$.
  Such edges could either join vertices in
  $S_1$ and   $S_4$, or in  $S_1$ and $S_5$, or in  $S_2$ and $S_5$,
  which require to ``jump'' over $\{S_2,S_3\}$ or  $\{S_3,S_4\}$ at least.
  We illustrate this in Figure   \ref{fig_upperLowerBound}.
\end{example}

The following theorem forms the basis for our lower bounds on the bandwidth.
\begin{thm}
  \label{prop:propBasicPartition}
  Let $A$ be an $n \times n $ adjacency matrix, and let
  $3 \leq k \leq n$ and $m \in \N^{k}$ be given with
  $\sum\limits_{i=1}^{k} m_{i} =n$. Let $1 \leq r \leq k-2$. If
  $$
  \min_{P \in \Pi_n}
  \frac{1}{2}
  \langle \transpose{P} A P,
  \overline{X} B_{r,k} \transpose{\overline{X}}  \rangle >0,
  ~\mbox{then}
  $$
  \begin{equation*}
    \bdw(A) > \min \{ m_{2} + \dotsm + m_{r+1}, ~
    m_{3} + \dotsm + m_{r+2}, \ldots, ~ m_{k-r} + \dotsm + m_{k-1}\}.
  \end{equation*}
  \begin{proof}
    If
    $\langle \transpose{P} A P,
    \overline{X} B_{r,k} \transpose{\overline{X}}  \rangle >0,$
    then some nonzero entry of $\transpose{P} A P$ is multiplied with
    a nonzero entry of  $\overline{X} B_{r,k} \transpose{\overline{X}}$.
    The nonzeros of this matrix closest to the main diagonal are in the
    positions
    $$
    (m_1, m_1 + \dotsm + m_{r+1} +1), \ldots,
    (m_1 + \dotsm +m_{k-r-1}, m_1 + \dotsm + m_{k-1} +1).
    $$
    As an illustration, these positions are marked with bullets in Figure   \ref{fig_upperLowerBound} below.
    The distances of these positions to the main diagonal are given by
    $$
    m_2 + \dotsm + m_{r+1}, \ldots, m_{k-r} + \dotsm + m_{k-1}.
    $$
    Therefore $\bdw(A)$ must be larger than the   smallest of these numbers.
  \end{proof}
\end{thm}
In case that the above minimum is zero, we have to consider
the zeros of $\overline{X} B_{r,k} \transpose{\overline{X}}$ with
largest possible distance to the main diagonal. These are marked with
crosses in Figure \ref{fig_upperLowerBound}.
\begin{thm}
  \label{prop:upperBound}
  Let $A$ be an $n \times n $ adjacency matrix, and let
  $3 \leq k \leq n$ and $m \in \N^{k}$ be given with
  $\sum\limits_{i=1}^{k} m_{i}
  =n$. Let $1 \leq r \leq k-2$. If
  $$\min_{X \in \mathcal{P}_{m}} \frac{1}{2} \langle A, XB_{r,k}\transpose{X} \rangle =0,
  ~\mbox{then}
  $$
  \begin{equation*}
    \bdw(A) < \max \{
    m_1 + m_2 + \dotsm + m_{r+1}, ~
    m_2 + m_3 + \dotsm + m_{r+2},
    \ldots, ~
    m_{k-r}+ \dotsm + m_{k}
    \}.
  \end{equation*}
\end{thm}
The proof is similar to Theorem \ref{prop:propBasicPartition} and is therefore
omitted.
In Figure 	\ref{fig_upperLowerBound}, we illustrate the
lower and upper bounds given by Theorems
\ref{prop:propBasicPartition} and \ref{prop:upperBound}, respectively.

\begin{figure}[!ht]
  \centering
  \includegraphics[scale=0.3]{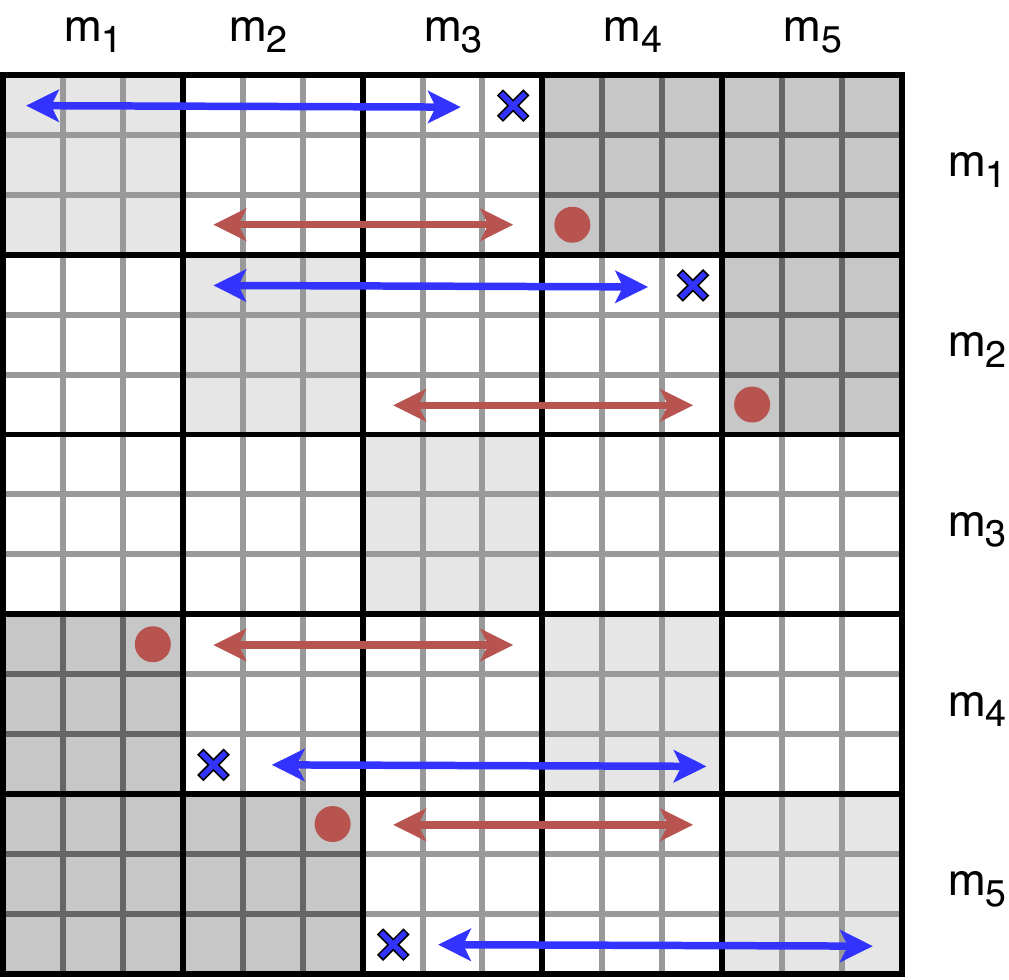}
  \caption{$A \in \cS_{15}$, $m=\transpose{(3,3,3,3,3)}$, and $r=2$. The crosses (bullets) indicate
  possible positions of the non-zero entries in terms of lower (upper) bounds. }
  \label{fig_upperLowerBound}
\end{figure}

The following   \textit{Minimal Partition Problem} ($\minpart$):
\begin{equation}\label{minPart}
  \minpart(m,r):=
  \min_{X \in \mathcal{P}_{m}} \frac{1}{2} \langle A, XB_{r,k}\transpose{X} \rangle
\end{equation}
serves as the basis to derive lower bounds on the bandwidth of $A$.
From a practical point of view, we are interested in selections of  $m$
where the minimum in Theorem \ref{prop:propBasicPartition} is attained
in each term. Some particular cases are summarized in the following
corollaries.
\begin{cor}	\label{cor:equiBdw}
  Let $A$ be an $n \times n $ adjacency matrix of $G$,  and let 	$3 \leq k \leq n$ and $m \in \N^{k}$ be given with
  $\sum_{i=1}^{k} m_{i} =n$. Let $r=1$. Further, suppose that $m_{2} = \dotsm = m_{k-1}$.

  If there exists $X \in P_{m}$ such that $\langle A, XB_{r,k}\transpose{X}
  \rangle >0$, then $\bdw(A) >  m_{2}$.
\end{cor}
\begin{cor} \label{cor:alternatingM}
  Let $A$ be an $n \times n $ adjacency matrix of $G$, and let $r=2$ and
  $m \in \N^{k}$ be given with  $\sum_{i=1}^{k} m_{i} =n$.
  Further, suppose  $m=\transpose{(m_1, m_2, m_3, m_2, \allowbreak m_3, \ldots, m_k)}$.

  If there exists $X \in \mathcal{P}_{m}$ such that $\langle A, XB_{r,k}\transpose{X} \rangle >0$,
  then $\bdw(A)> m_{2}+m_{3}.$
\end{cor}
By cyclically repeating the sizes, we can insure that the minimum in
Theorem  \ref{prop:propBasicPartition} is attained in each term
simultaneously as above also for values $r>2$.

\subsection{Relation to Prior Work}
We present below two important special cases of our new modelling approach and their relation to  prior work.

\paragraph{The case $k=3$}
Given $k=3$ the only allowable choice for $r$ is $r=1$ and therefore
the only nonzero elements in $B_{{1,3}}$ are $b_{1,3}=b_{3,1}=1$.
Hence for $m=\transpose{(m_{1},~ m_{2},~m_{3})}$
Theorem \ref{prop:propBasicPartition} states that
if there exists $X \in \mathcal{P}_{m}$ such that  $ \langle A,XB_{1,3}\transpose{X} \rangle>0$,  then
$\bdw(A)>m_{2}$.
This observation is used in \cite{HelReMoPo:95} to derive lower bounds on $\bdw(A)$, and is further refined in \cite{PoRe:07,DamSotirov15}.

\paragraph{The case $k=n$}
Another notable case occurs for $k=n$, which implies that  $m_{1}= \dotsm = m_{n} =1$.
Hence,  for any $r \in \{1, \ldots, n-2 \}$ it follows from
Theorem \ref{prop:propBasicPartition}   that
$\bdw(A)>r$, if there exists a partition matrix $X \in \mathcal{P}_{m}$ such that $\langle A, XB_{r,n}\transpose{X}  \rangle >0$.
However, in this case the basic partition matrix becomes the identity matrix of
rank $n$, i.e., $\overline{X}=I_n$. Thus, $X$ becomes a permutation matrix  $Q\in \Pi_n$
and we recover the statement
$$
\min\limits_{Q \in \Pi_n}  \langle \transpose{Q} A Q, B_{r,n} \rangle
>	0 \Rightarrow  \bdw(A) > r,
$$
from \eqref{eq:qap}.
This approach is used e.g., in \cite{dKlSoNa} to derive lower
bounds on $\bdw(A)$.

In Summary, we have shown that once the $\minpart$ problem has a
positive value for given $B_{r,k}$ and $m$, we get a nontrivial
lower bound on the
bandwidth from Theorem \ref{prop:propBasicPartition}.
The $\minpart$ problem
is itself NP-complete, so our strategy is to consider tractable lower bounds
for the $\minpart$ problem.
If some lower bound turns out to be positive for
given $r$ and $m$, then clearly $\minpart$ has a positive value, and
our bounding argument can be applied. In the following section we consider
relaxations of $\minpart$, based on semidefinite optimization.

\section{SDP models} \label{sect:SDPmodels}

In this section, we derive several Semidefinite Programming relaxations for
the Minimal Partition Problem.
Our first two SDP relaxations are obtained by matrix lifting and therefore  have
matrix variables of order ${\mathcal O}(n\cdot k)$, while the third relaxation has $k$  matrix variables of order $n$.

\subsection{SDP model in  ${\mathcal S}_{n\cdot k+1}$}

In this section, we derive an SDP relaxation whose matrix variable is of order $n\cdot k+1$.

Let $X \in \mathcal{P}_m$ be a partition matrix, see \eqref{PM}.
Let $x_1,\ldots,x_k$ be the columns of $X$,
i.e., $X=\begin{bmatrix}
x_1 & \cdots & x_k
\end{bmatrix}$,
and
$x:=\vect(X)\in \mathbb{R}^{n \cdot k}$.
Now, the constraint $Y =x \transpose{x}$ may be weakened to $Y -x \transpose{x} \succeq 0$
which is well-known to be equivalent to the following convex constraint
\[
Z:=\begin{bmatrix}
Y & x\\
\transpose{x} & 1
\end{bmatrix} \succeq 0.
\]
Further, we  use the following block notation for $Z \in \mathcal{S}_{n\cdot k +1}$:
\[
Z= \begin{bmatrix}
X_{1} 													& X_{12} 														&  \ldots 						&  X_{1k} 						& x_1  \\[1ex]
\transpose{X_{12}} 							& X_{2} 														&   \ldots 						&  X_{2k}  						& x_2  \\[1ex]
\vdots													&			\vdots 												& 					 					& 		\vdots				  & 		\vdots \\[1ex]
\transpose{X_{1k}}              &  \transpose{X_{2k}} 	             &  \ldots            & X_{k}		            & x_k  \\[1ex]
\transpose{x_{1}} 							& \transpose{x_{2}}									&  \ldots 						&  \transpose{x_{k}}	& 1
\end{bmatrix},
\]
where $X_i$ corresponds to $x_i \transpose{x_i} ,~  i=1,\ldots, k$, and $X_{ij}$ to $x_i \transpose{x_j}, ~  i\neq j,~  i,j=1,\ldots,k$.

For any $X_i, ~ i=1,\ldots,k$, we have  $\diag(X_i)=\diag(x_i \transpose{x_i})=x_i$ and thus $\trace\left( X_i \right)= \transpose{x_i}e_n  = m_i$.
For all   $i=1,\ldots,k$  we have:
\[
\langle J_n, X_i \rangle
= \trace \big( e_n \transpose{e_n} x_i \transpose{x_i} \big)
= \trace \big( (\transpose{x_i}e_n)^2 \big)
= {m_i}^2.
\]
Similarly, we have
\begin{align*}
  \langle J_n, X_{ij} + \transpose{X_{ij} } \rangle
  &= \trace \big( J_n X_{ij} + J_n \transpose{X_{ij}} \big)\\
  &= 2 \cdot  \trace \big( e_n \transpose{e_n} x_i \transpose{x_j} \big)
  = 2  m_i m_j, \quad \forall i,j.
\end{align*}
From orthogonality of vectors $x_i$, $i=1,\ldots,k$, it follows $\diag \big( X_{ij} \big)=0$.

Let us  describe the matrix \eqref{Brk} as the sum of symmetric matrices having only two non-zero entries, i.e.,
$B_{r,k}=  \sum\limits_{|u-v| >r} \left(
e_{k}^{u} \transpose{{e_{k}^v}}  + e_{k}^{v} \transpose{{e_{k}^u}} \right).$
Hence, we derive
\begin{align*}
  X B_{r,k} \transpose{X}
  =&
  \sum\limits_{|u-v| >r}
  \left(
  X  e_{k}^{u} \transpose{{e_{k}^{v}}} \transpose{X}  +  X e_{k}^{v} \transpose{{e_{k}^{u}}} \transpose{X}  \right)\\
  =&
  \sum\limits_{|u-v| >r} \left(
  x_{u} \transpose{x_{v}}  + x_{v} \transpose{x_{u}} \right) =
  \sum\limits_{|u-v| >r} \left( X_{uv} +X_{vu} \right).
\end{align*}
Therefore, we can rewrite the Minimal Partition Problem, see \eqref{minPart}, as:
\begin{align*}
  \min\limits_{X \in \mathcal{P}_m}
  \frac{1}{2} \langle  A, X B_{r,k} \transpose{X}  \rangle
  &=
  \min\limits_{X \in \mathcal{P}_m}  \frac{1}{2} \langle A, \sum\limits_{|u-v| >r} X_{uv} +X_{vu} \rangle
  =
  \min\limits_{X \in \mathcal{P}_m}  \sum\limits_{|u-v| >r} \langle A, X_{uv}     \rangle
  .
\end{align*}

Finally, we collect all above mentioned constraints  and propose the following model for the Minimal Partition Problem
based on the matrix lifting approach.
\begin{subequations}
  \begin{align}
    \min \quad & \sum\limits_{|u-v| >r} \langle A, X_{uv}     \rangle,
    \\
    \text{s.t.} \quad
    &	\diag\left(X_i\right)=x_i, \quad i=1,\ldots,k,
    \label{cond:ml1diag}	\\
    & \diag\left(X_{ij}\right) = 0 ,
    \quad i \neq j,~ i,j=1,\ldots k,
    \label{cond:ml1diag2}
    \\
    &	\trace\left(X_i\right) = m_i  , \quad i=1,\ldots, k,
    \label{cond:ml1sum1} \\
    &	\langle J_n, X_i \rangle = {m_i}^2, \quad i=1,\ldots, k,
    \label{cond:ml1sum3} \\
    &	\langle J_n, X_{ij}  + \transpose {X_{ij}} \rangle  = 2 m_i  m_j ,
    \quad i \neq j,
    ~ i,j=1,\ldots, k,
    \label{cond:ml1mm} \\
    & Z= \begin{bmatrix}
    X_{1} 													& X_{12} 														&  \ldots 						&  X_{1k} 						& x_1  \\[1ex]
    {X_{21}} 							& X_{2} 														&   \ldots 						&  X_{2k}  						& x_2  \\[1ex]
    \vdots													&			\vdots 												& 					 					& 		\vdots				  & 		\vdots \\[1ex]
    {X_{k1}}              &  {X_{k2}} 	             &  \ldots            & X_{k}		            & x_k  \\[1ex]
    \transpose{x_{1}} 							& \transpose{x_{2}}									&  \ldots 						&  \transpose{x_{k}}	& 1
  \end{bmatrix}
  \succeq 0.
  \label{cond:ml1symmetry}
\end{align}
\label{mod:matrixliftinglarge}
\end{subequations}
Here  $Z \in \mathcal{S}_{k n +1}^{+}$.
The feasible region of the above SDP relaxation equals  the feasible  region  of the SDP relaxation for the graph partition problem derived by Wolkowicz and Zhao \cite{WOLKOWICZ1999461}.
In order to further improve the relaxation, one can add nonnegativity constraints.

Below, we analyze the feasible region of the model \eqref{mod:matrixliftinglarge}.

\begin{lem} \label{nullspace}
  Let $Z$ satisfy  \eqref{cond:ml1diag}, \eqref{cond:ml1diag2}, \eqref{cond:ml1sum1},  \eqref{cond:ml1sum3}, and \eqref{cond:ml1symmetry}. Then
  $$ \underbrace{
  \begin{pmatrix}
    e_n \\ 0_n \\ \vdots \\ 0_n \\ -m_1
    \end{pmatrix},
    \begin{pmatrix}
      0_n \\ e_n \\ \vdots \\ 0_n \\ -m_2
      \end{pmatrix},
      \ldots,
      \begin{pmatrix}
        0_n \\ 0_n \\ \vdots \\ e_n \\ -m_k
      \end{pmatrix}
      }
      _{k ~\text{vectors}}
      ,
      \underbrace{
      \begin{pmatrix}
        I_n \\ I_n \\ \vdots \\ I_n \\ -\transpose{e_n}
      \end{pmatrix}
      }
      _{n ~ \text{vectors}}
      $$
      spans the nullspace of $Z$.
    \end{lem}
    For a proof we refer the reader to
    \cite[Lemma 10 and Section 5.2]{Rendl2017}   as well as to
    \cite{WOLKOWICZ1999461}.
    We observe in particular  that this result holds independent of
    \eqref{cond:ml1mm}.

    Note that the vectors from Lemma \ref{nullspace} correspond to a $(n\cdot k +1) \times (n+k)$ matrix.
    As the sum of the first $k$ columns is equal to the sum of the last $n$ columns, the nullspace of $Z$ has dimension $n+k-1$.

    \begin{lem} \label{lem:matrixliftinlemma2}
      Let $Z$  satisfy \eqref{cond:ml1diag}, \eqref{cond:ml1diag2}, \eqref{cond:ml1sum1},  \eqref{cond:ml1sum3}, and \eqref{cond:ml1symmetry}.
      Then
      \[
      \left\{
      \begin{array}{lllll}
        X_1 &+ X_{12} &+ \dotsm &+ X_{1k} &= x_1 \transpose{e_n} \\
        \vdots &      &        &           &~  \vdots\\
        X_{k1} &+ X_{k1} &+ \dotsm &+ X_{k} &= x_k \transpose{e_n} \\
        x_1 &+ x_2 &+ \dotsm &+ x_k &= e_n
      \end{array}
      \right.
      .
      \]
    \end{lem}
    Again, we refer the reader to \cite[Section 5.2]{Rendl2017}, and  \cite{WOLKOWICZ1999461} for a formal proof.
    As a consequence of the previous lemma, the block
    $\begin{bmatrix}
    X_{k1} & X_{k2} & \ldots& X_ {k,k-1}& X_k & x_k
  \end{bmatrix}
  $
  is determined by
  $X_1,\ldots, ~X_{k-1}$,
  $X_{ij}$, ($i\neq j, ~i,j=1,\ldots,k-1$),
  and
  $x_1,\ldots,x_{k-1}$.
  Hence, matrix $Z$ can be reduced by one block of rows and their corresponding columns without loss of information.
  This leads us to the reduced SDP model presented in the following section.

  One can also derive  the Slater feasible version
  of the SDP relaxation \eqref{mod:matrixliftinglarge} by exploiting a basis of the orthogonal complement to the nullspace of $Z$ given in  Lemma \ref{nullspace}.
  For details see e.g., \cite{ZhKaReWo:98,Rendl2017}.
  The Slater  feasible version may be efficiently solved by using the Alternating Direction Method of Multipliers (ADMM) as described in \cite{HenryADMM}.
  The ADMM is  a first-order method for convex problems that decomposes an optimization problem into subproblems that may be easier to solve.

  \subsection{Reduced SDP Model in $\mathcal{S}_{n\cdot (k-1)+1}$}

  In this section, we provide an SDP relaxation that is equivalent to the one from the previous subsection, but contains less variables.
  In particular, based on Lemma  \ref{lem:matrixliftinlemma2}, we propose the following SDP relaxation for the Minimal Partition Problem.
  \begin{subequations}
    \begin{align}
      \min \quad & \sum\limits_{|u-v| >r} \langle A, X_{uv}     \rangle,
      \\
      \text{s.t.} \quad
      &	\diag\left(X_i\right)=x_i, \quad i=1,\ldots,k-1,
      \label{cond:ml2diag}	\\
      & \diag\left(X_{ij}\right) = 0 ,
      \quad i \neq j,~ i,j=1,\ldots k-1,
      \label{cond:ml2diag2}
      \\
      &	\trace\left(X_i\right) = m_i  , \quad i=1,\ldots, k-1,
      \label{cond:ml2sum1} \\
      &	\langle J_n, X_i \rangle = {m_i}^2, \quad i=1,\ldots, k-1,
      \label{cond:ml2sum3} \\
      &	\langle J_n, X_{ij}  + \transpose {X_{ij}} \rangle  = 2 m_i  m_j ,
      \quad i \neq j,
      ~ i,j=1,\ldots, k-1, \label{cond:ml2mm} \\
      & \widetilde{Z}= \begin{bmatrix}
      X_{1} 													& X_{1,2} 														&  \ldots 						&  X_{1,k-1} 						& x_1  \\[1ex]
      \transpose{X_{1,2}} 							& X_{2} 														&   \ldots 						&  X_{2,k-1}  						& x_2  \\[1ex]
      \vdots													&			\vdots 												& 					 					& 		\vdots				  & 		\vdots  \\[1ex]
      \transpose{X_{1,k-1}}              &  \transpose{X_{2,k-1}} 	             &  \ldots            & X_{k-1}		            & x_{k-1}   \\[1ex]
      \transpose{x_{1}} 							& \transpose{x_{2}}									&  \ldots 						&  \transpose{x_{k-1}}	& 1
    \end{bmatrix}
    \succeq 0.
    \label{cond:ml2symmetry}
  \end{align}
  \label{mod:ml2}
\end{subequations}
\noindent
Here $\widetilde{Z} \in \mathcal{S}_{n \cdot(k-1)+1}^{+}$.
Note that the nullspace of the reduced matrix $\widetilde{Z}$ has rank $k-1$. We show below that the SDP relaxation \eqref{mod:ml2}
is equivalent to \eqref{mod:matrixliftinglarge}.
{
The number of equations in this SDP is still
${\mathcal O}(n \cdot k)$, but we saved about $n$ equations as
compared to the original model.
}

Additional sign constraints
\begin{equation}
  \label{eq:SDPnonneg}
  X_{uv} \geq 0, \quad |u-v|> r
\end{equation}
insure that the lower bound from this model is always nonnegative.

\begin{lem}
  From \eqref{cond:ml2diag} --  \eqref{cond:ml2symmetry}
  follow
  \eqref{cond:ml1diag} -- \eqref{cond:ml1symmetry}.
  \begin{proof}
    \proofstep{1}      From  Lemma  \ref{lem:matrixliftinlemma2}  directly
    follows that, given $\widetilde{Z}$,  the ``missing'' entries of $Z$  can be expressed by:
    \begin{align}
      x_k &= e_n -x_1- \cdots - x_{k-1} \geq 0, \nonumber\\
      X_{ik} &= x_i \transpose{e_n} - X_i - \sum\limits_{\substack{j=1 \\ i\neq j}}^{k-1} X_{ij}
      , \quad i=1,\ldots,k-1   \nonumber,\\
      X_{k} &= x_k \transpose{e_n}  - \sum\limits_{j=1}^{k-1} X_{kj}.  \nonumber
    \end{align}
    Nonnegativity of $x_k$ follows from \eqref{cond:ml2diag2} and \eqref{cond:ml2symmetry}.

    \noindent
    \proofstep{2}

    \paragraph{Constraint \eqref{cond:ml1symmetry}}
    From \cite[Section 5]{Rendl2017}, we know that under
    \eqref{cond:ml2diag} -- \eqref{cond:ml2symmetry}
    it holds that
    $$
    \widetilde{Z} \succeq 0 ~~ \wedge ~~  Z=\overline{W} U \transpose{\overline{W}}  \Rightarrow Z\succeq 0,
    $$
    where
    \begin{equation}\label{Wb}
      \overline{W}:=
      \begin{bmatrix}
        e_n    &  0_n    & \cdots  &  0_n      & I_n   \\
        0_n    &  e_n    &  \cdots   &  0_n    & I_n   \\
        \vdots &  \vdots  &  \cdots  & \vdots  &\vdots \\
        0_n    &   0_n   &   \cdots   & e_n    & I_n   \\
        -m_1   &  -m_2  &   \cdots    & -m_k   & -e_n  \\
        \end{bmatrix}^{\bot}.
      \end{equation}
      Hence, it holds \eqref{cond:ml1symmetry}.

      \paragraph{Constraint \eqref{cond:ml1diag}}
      In addition to \eqref{cond:ml2diag}, $\diag(X_k)= x_k$ must hold. In particular, from {\bf Step 1} it follows
      \begin{align*}
        \diag\big(X_k \big) {=}
        \diag \big( x_k \transpose{e_n} -  \sum\limits_{j=1}^{k-1} X_{kj} \big)
        =  x_k -  \sum\limits_{j=1}^{k-1}  \diag (X_{kj})
        = x_k.
      \end{align*}

      \paragraph{Constraint \eqref{cond:ml1diag2}}
      In addition to \eqref{cond:ml2diag2}, $\diag(X_{ik})=0,$ $i=1,\ldots,k-1$, must hold. Again, by using {\bf Step 1} we have:
      \[
      \diag \big( X_{ik}\big) = \diag\big( x_i \transpose{e_n}  - X_i - \sum\limits_{\substack{j=1 \\ i\neq j}}^{k-1}  X_{ij}\big)=0.
      \]
      \paragraph{Constraint  \eqref{cond:ml1sum1}}
      From  \eqref{cond:ml2diag} and  {\bf Step 1}  we have \eqref{cond:ml1diag}. \\Thus, from $\diag(X_k)=x_k$ it follows $\trace \big(X_k\big)=m_k$.

      \paragraph{Constraint  \eqref{cond:ml1sum3}}
      From  $\langle  J_n, x_k\transpose{e_n} \rangle =  m_k \cdot n$
      and $\langle J_n ,  X_{kj} \rangle=   \langle e_n \transpose{e_n}, x_k \transpose{x_j} \rangle =m_{j}\cdot m_{k}$, we have
      \[
      \langle J_n, X_k \rangle =   \langle J_n, x_k \transpose{e_n} - \sum\limits_{j=1}^{k-1} X_{kj} \rangle =
      m_k \big( n- \sum\limits_{j=1}^{k-1} m_j\big)= m_k^{2}.
      \]

      \paragraph{Constraint  \eqref{cond:ml1mm}}
      In addition to   \eqref{cond:ml2mm}, $\langle J_n, X_{ik}+\transpose{X_{ik}} \rangle= 2 m_{i} m_{k}, ~ i=1,\ldots,k-1$,
      must hold.
      \begin{align*}
        \langle J_n, X_{ik}+ \transpose{X_{ik}} \rangle  &= 2 \cdot  \langle J_n, X_{ik} \rangle =
        ~2 \cdot \big[
        \langle J_n, x_i \transpose{e_n} \rangle - \langle J_n, X_i \rangle
        -\sum\limits_{\substack{j=1 \\ i\neq j}}^{k-1}  \langle J_n, X_{ij}  \rangle
        \big]\\
        &=~2 m_{i} m_{k}.
      \end{align*}
    \end{proof}
  \end{lem}

  Note that the inverse to the one in the lemma follows directly.
  To make the SDP relaxation \eqref{mod:ml2} with additional nonnegativity constraints
  equivalent to  SDP relaxation \eqref{mod:matrixliftinglarge} with additional nonnegativity constraints,
  we need to add nonnegativity constraints to the ``missing'' blocks
  $[X_{k1} ~ X_{k2} ~ \dotsm ~ X_ {k,k-1} ~~ X_k $ $~~ x_k]$ in \eqref{mod:ml2}.
  In particular, we have the following proposition.
  \begin{prop} \label{addConstr}
    The SDP relaxation \eqref{mod:matrixliftinglarge} with additional  constraints ${Z}\geq 0$
    is equivalent to the SDP relaxation   \eqref{mod:ml2}  with additional  constraints $\widetilde{Z}\geq 0$  and
    \begin{align*}
      & 1 -\sum_{r=1}^{k-1} (X_r)_{i,i} - \sum_{r=1}^{k-1} (X_r)_{j,j} + \sum_{r=1}^{k-1}\sum_{p=1}^{k-1} (X_{rp})_{i,j} \geq 0, ~~i>j
      ,  \\
      & (X_r)_{i,i} - \sum_{l=1}^{k-1} (X_{lr})_{i,j} \geq 0, \quad i\neq j,~ r\in \{1, \ldots, k-1\},
    \end{align*}
    where $ i,j=1,\ldots,n$.
  \end{prop}

  In Section \ref{sect:numerResults}, we demonstrate the strength of our SDP relaxation.

  \section{Computational Experiments}
  \label{sect:numerResults}

  \subsection{Solving the SDP relaxation}

  The partition-based lower bounds for the bandwidth problem lead to
  semidefinite programs with  one  matrix of dimension
  $n \cdot(k-1) + 1$, see \eqref{mod:ml2}.
  The resulting relaxations can be solved using
  standard SDP packages such as SDPT3 only for limited values of $n$ and $k$.
  % We now focus on the strongest model which has one matrix variable of order $n \cdot(k-1)+1$ and
  % roughly $n \cdot k^{2}$ equality constraints.
  We also consider nonnegativity
  constraints which add another $O(n^{2}k^{2})$ potentially violated sign
  constraints to our relaxation. Interior-point based methods for such a
  scenario turn out to
  be too slow. Hence, we propose to use the
  ADMM method, which works well for SDPs with simple sign constraints.
  To use the ADMM, we use the Slater feasible version of the SDP
  relaxation \eqref{mod:matrixliftinglarge} as described in the
  previous section.
  The resulting SDP relaxation has a matrix variable of order $(k-1)\cdot (n-1)+1$,
  see e.g., \cite{WOLKOWICZ1999461}.
  Then, we proceed in the same  manner as described in   \cite{HenryADMM,HaoSotirov:19}.

  \subsection{Strength of the partition bounds}

  As a first experiment we investigate the quality of the
  SDP relaxations   \eqref{mod:matrixliftinglarge}
  and \eqref{mod:ml2}
  to assess
  $$
  \mbox{minPart}(m,r) > 0
  $$
  for given $m$ and $r$. We recall that minPart($m,r$) denotes the number of
  edges in the minimal partition specified by $m$ and $r$, see
  \eqref{minPart}.
  We are primarily interested in parameter settings for $m$ and $r$
  where minPart$(m,r) > 0$ but small. For such values of $m$ and $r$
  it is a nontrivial task to prove positive lower bounds for
  minPart using our SDP models.

  \subsubsection{Test problems}

  We investigate the practical performance of our lower bounds on the
  following classes of graphs.

  \paragraph{Torus graphs}
  For given integer $k$ the  torus graph $T_{k}$ has $k^{2}$ vertices which we
  label by $(i,j)$ for $i,j \in \{1,
  \ldots,k\}$. We introduce  ``vertical'' edges of the form $[ (i,j), (i+1,j)]$
  for $1\leq i \leq k-1$ and $[(1,j),(k,j)]$. Altogether there
  are $k^{2}$ such edges. In a similar way we add ``horizontal'' edges of the
  form $[(i,j),(i,j+1)]$ for $j<k$ together with
  $[(i,1),(i,k)]$. This graph therefore has
  $n:=k^{2}$ vertices and $2n$ edges.
  These graphs are interesting for the following reason.
  They are extremely sparse ($n$ vertices and $2n$ edges), but
  their bandwidth is quite large. Namely,
  it is known that $\bdw(T_{k})=2k-1$, see e.g., \cite{BALOGH200643,LiTapShen}.

  \paragraph{Torus graphs plus Hamiltonian path}
  Here we start out with the torus graph $T_{k}$, choose a labeling of
  its vertices yielding a bandwidth of size $2k$, and add the
  Hamiltonian path from the first to the last vertex in this labeling.
  The resulting graph is denoted by $TH_{k}$. It is still sparse having
  roughly $3|V(TH_{k})|$ edges and bandwidth again at most $2k$.

  \paragraph{Hypercubes }
  The Hamming graph $H(d,q)$ is the Cartesian product of $d$ copies of the
  complete graph $K_q$. The Hamming graph
  $H(d,2)$ is also known as the hypercube (graph) $Q_d$. Thus, the hypercube
  graph  $Q_d$ has $2^d$ vertices.
  The bandwidth of the hypercube graph was determined
  by Harper \cite{Harper66}
  and is given by the following expression:
  \[
  \bdw(Q_d) = \sum\limits_{i=0}^{d-1} \binom{i}{ \lfloor \frac{i}{2} \rfloor}.
  \]
  We use the  hypercube graphs $Q_d$ to test the quality of our
  partition bounds.

  \subsubsection{Computations}

  \paragraph{Torus graphs}
  In the tables to follow we always provide the following information.
  The first block of data contains the vector $m$ of cardinalities for the partition blocks.
  We consider partitions into $k \in \{4,5,6 \}$ blocks. We set $r=1$ and ask that
  $
  m_2 = m_3 = \dotsm = m_{k-1}.
  $

\begin{figure}
    \includegraphics[width=0.3\textwidth]{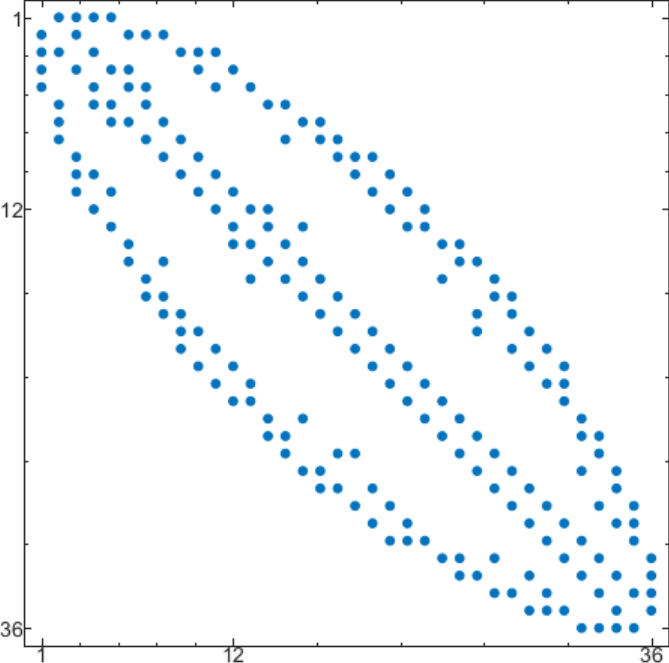}
    \hfill
    \includegraphics[width=0.3\textwidth]{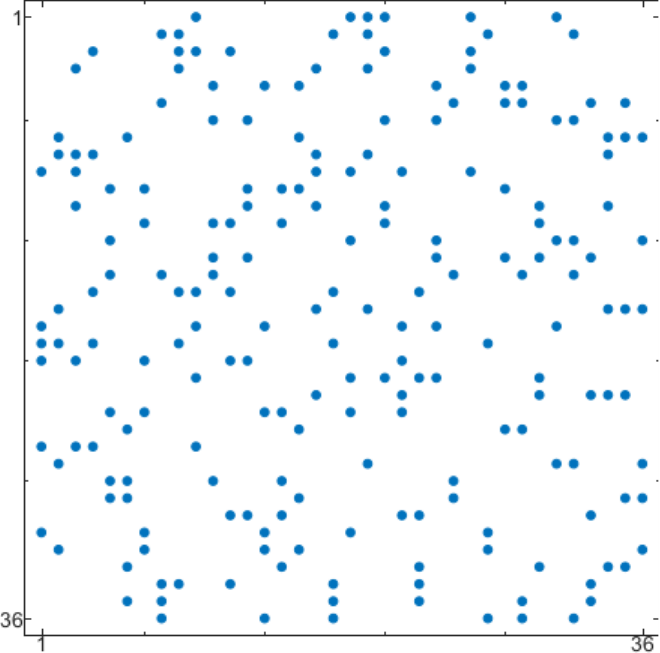}
    \hfill
    \includegraphics[width=0.3\textwidth]{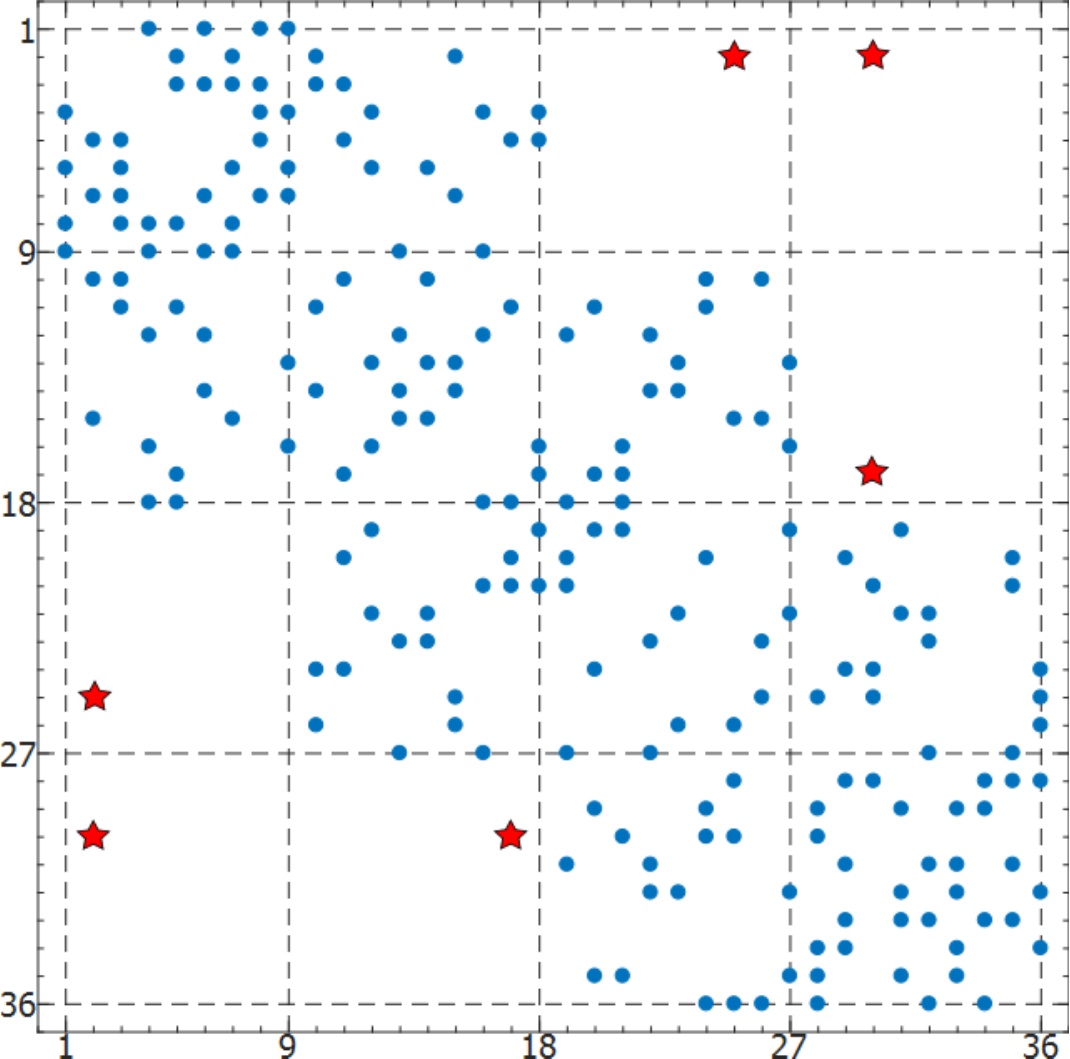}

    \caption{Illustration of the $TH_6$ graph. On the left, we show the unpermuted graph, in the center,
    the permuted graph is shown, on the right, the obtained solution of the $\minpart$ problem with
    $m=\transpose{(9,9,9,9)}$ is shown. The value of $\minpart$ is $3$, the corresponding entries are indicated by stars. }
    \label{fig:stages}
  \end{figure}

  The sizes $m_1$ and $m_k$ are chosen such that $\sum_{i=1}^{k} m_i =n$ and
  $| m_1 -m_k | \leq 1.$ Next we provide upper and lower bounds for the
  Minimal Partition Problem.  The
  upper bound (ub)
  is obtained by running a
  standard Simulated Annealing heuristic \cite{BURKARD1984169} to find a
  good partition.
  The lower bound (lb)
  is obtained from the     SDP relaxation   \eqref{mod:matrixliftinglarge}
  with all nonnegativity constraints included.
  Our main interest lies in values of $m$, where the obtained lower bound is
  nontrivial, i.e., $lb>0$.
  We give an illustration of the obtained solutions in Figure \ref{fig:stages}.

  First, we consider Table \ref{tab:T7}, which contains computational results
  for the Torus graph $T_7$.
  Initially, we consider 4 blocks with $m_2= m_3=8$ leading to a
  lower bound $lb>1.23$. Hence, Corollary \ref{cor:equiBdw} allows us to
  conclude that $\bdw(T_7)>8$.
  We next try $m_2=m_3=9$ where we only obtain the trivial lower bound of $0$.
  Therefore,
  we get no further restriction on $\bdw(T_7)$ from 4-partitions.
  The 5-partition with $m_2 = m_3 = m_4= 9 $ however yields a positive lower
  bound and therefore $\bdw(T_7)>9$.
  Also, 6-partitions, given in the last block of the Table \ref{tab:T7},
  do not lead to a further tightening of $\bdw(T_7)$.

  \begin{table}[!ht]
    \begin{center}
      \begin{tabular}{rrrrrr|r|r}
        \multicolumn{8}{c}{$T_{7}  \qquad (n=49)$} \bigstrut\\
        \hline\hline
        $m_1$ & $m_2$ & $m_3$ &$m_4$ & $m_5$ & $m_6$ & $ub $ & $lb$ \\
        \hline\hline
        16 & 8 & 8 & 17 &   &    & 6 & 1.23 \\
        15 & 9 & 9 & 16 &   &    & 5 & -    \\ \hline
        11 & 9 & 9 &  9 & 11 &   & 6 & 0.68 \\
        9 & 10& 10& 10 & 10 &   & 5 & - \\ \hline
        6 & 9 & 9 &  9 &  9 & 7 & 6 & 0.56 \\
        4 & 10& 10& 10 & 10 & 5 & 4 & - \\ \hline
      \end{tabular}
      \caption{Torus graph $T_{7}$.}
      \label{tab:T7}
    \end{center}
  \end{table}

  The results for the Torus graphs $T_8, T_9$, and $T_{10}$ are summarized in
  Table \ref{tab:t8910}.
  We proceed as before and consider partitions with $k \in \{4, 5, 6 \}$.
  We can prove a lower bound of $11$ for $\bdw(T_8)$ and $\bdw(T_9)$.
  It turns out that proving positive lower bounds for our partition problems gets
  increasingly difficult as either $n$ or $k$ increases.
  For  $T_ {10}$, the use of a 6-partition allows us to prove a lower bound of $14$.

  \begin{table}[!ht]
    \begin{center}
      \begin{tabular}{rrrrrr|r|r}
        \multicolumn{8}{c}{$T_{8}  \qquad (n=64)$} \bigstrut\\
        \hline\hline
        $m_1$ & $m_2$ & $m_3$ &$m_4$ & $m_5$ & $m_6$ & $ub $ & $lb$ \\
        \hline\hline
        23 & 9 & 9 & 23 &   &    & 7 & 1.01 \\
        22 &10 &10 & 22 &   &    & 6 & -    \\ \hline
        17 &10 &10 & 10 & 17 &   & 7 & 0.84 \\
        15 & 11& 11& 11 & 16 &   & 7 & - \\ \hline
        12 &10 &10 & 10 & 10 &12 & 8 & 0.99 \\
        10 & 11& 11& 11 & 11 &10 & 6 & - \\ \hline
        \multicolumn{8}{c}{$T_{9}  \qquad (n=81)$} \bigstrut\\
        \hline\hline
        $m_1$ & $m_2$ & $m_3$ &$m_4$ & $m_5$ & $m_6$ & $ub $ & $lb$ \\
        \hline\hline
        31 & 9 & 9 & 32 &   &    & 9 & 1.53 \\
        30 &10 &10 & 31 &   &    & 8 & -    \\ \hline
        25 &10 &10 & 10 & 26 &   &10 & 1.63 \\
        24 & 11& 11& 11 & 24 &   & 9 & - \\ \hline
        20 &10 &10 & 10 & 10 &21 & 9 & 1.91 \\
        18 & 11& 11& 11 & 11 &19 & 9 & - \\ \hline
        \multicolumn{8}{c}{$T_{10}  \qquad (n=100)$} \bigstrut\\
        \hline\hline
        $m_1$ & $m_2$ & $m_3$ &$m_4$ & $m_5$ & $m_6$ & $ub $ & $lb$ \\
        \hline\hline
        41 & 9 & 9 & 41 &   &    &11 & 1.62 \\
        40 &10 &10 & 40 &   &    &10 & -    \\ \hline
        32 &12 &12 & 12 & 32 &   &10 & 0.68 \\
        30 & 13& 13& 13 & 31 &   & 9 & - \\ \hline
        24 &13 &13 & 13 & 13 &24 &10 & 0.52 \\
        22 & 14& 14& 14 & 14 &22 &10 & - \\ \hline
      \end{tabular}
      \caption{Torus graphs $T_8,~ T_9,~ T_{10}$.}
      \label{tab:t8910}
    \end{center}
  \end{table}

  \medskip\medskip

  As a second experiment, we consider the graphs $TH_7, \ldots, TH_{10}$
  consisting of the union of the Torus graph and a Hamiltonian path such that
  $\bdw(TH_k) \leq 2k$ is insured.  The results are summarized in Table \ref{tab:torusHam}.
  Compared to the Torus graphs we get slightly stronger lower bounds even though
  these graphs are still quite sparse, with $|E(TH_k)| < 3|V(TH_k)|$.
  Again, we see increasing gaps between lower and upper bounds as
  the number of nodes of the graph increases.

  \begin{table}[!ht]
    \begin{center}
      \begin{tabular}{rrrrrr|r|r}
        \multicolumn{8}{c}{$TH_{7}  \qquad (n=49)$} \bigstrut\\
        \hline\hline
        $m_1$ & $m_2$ & $m_3$ &$m_4$ & $m_5$ & $m_6$ & $ub $ & $lb$ \\
        \hline\hline
        14 &10 &10 & 15 &   &    & 5 & 0.87 \\
        13 &11 &11 & 14 &   &    & 3 & -    \\ \hline
        8 &11 &11 & 11 &  8 &   & 2 & 0.18 \\
        6 & 12& 12& 12 &  7 &   & 1 & - \\ \hline
        2 &11 &11 & 11 & 11 & 3 & 3 & 0.07 \\
        \hline\hline
        \multicolumn{8}{c}{$TH_{8}  \qquad (n=64)$} \bigstrut\\
        \hline\hline
        $m_1$ & $m_2$ & $m_3$ &$m_4$ & $m_5$ & $m_6$ & $ub $ & $lb$ \\
        \hline\hline
        21 &11 &11 & 21 &   &    & 7 & 0.76 \\
        20 &12 &12 & 20 &   &    & 5 & -    \\ \hline
        14 &12 &12 & 12 & 14 &   & 6 & 0.64 \\
        12 &13 &13 & 13 & 13 &   & 3 & - \\ \hline
        8 & 12& 12& 12 & 12 & 8 & 6 & 0.35 \\
        6 &13 &13 & 13 & 13 & 6 & 3 & - \\
        \hline\hline
        \multicolumn{8}{c}{$TH_{9}  \qquad (n=81)$} \bigstrut\\
        \hline\hline
        $m_1$ & $m_2$ & $m_3$ &$m_4$ & $m_5$ & $m_6$ & $ub $ & $lb$ \\
        \hline\hline
        28 &12 &12 & 29 &   &    &10 & 0.96    \\
        27 &13 &13 & 28 &   &    & 7 & -    \\ \hline
        21 &13 &13 & 13 &21 &    & 8 & 1.12    \\
        19 &14 &14 & 14 &20 &    & 6 & -    \\ \hline
        14 &13 &13 & 13 &13 &15  &10 & 1.34    \\
        12 &14 &14 & 14 &14 &13  & 7 & -    \\
        \hline\hline
        \multicolumn{8}{c}{$TH_{10}  \qquad (n=100)$} \bigstrut\\
        \hline\hline
        $m_1$ & $m_2$ & $m_3$ &$m_4$ & $m_5$ & $m_6$ & $ub $ & $lb$ \\
        \hline\hline
        37 &13 &13 & 37 &   &    &11 & 0.64    \\
        36 &14 &14 & 36 &   &    & 9 & -    \\ \hline
        29 &14 &14 & 14 &29 &    &11 & 1.20    \\
        27 &15 &15 & 15 &28 &    & 9 & -    \\ \hline
        22 &14 &14 & 14 &14 &22  &11 & 1.64    \\
        20 &15 &15 & 15 &15 &20  & 9 & -    \\ \hline
      \end{tabular}
      \caption{Torus graphs plus Hamiltonian paths $TH_{7},~TH_{8}, ~TH_{9}, ~TH_{10}$.}
      \label{tab:torusHam}
    \end{center}
  \end{table}

  \begin{table}[!ht]
    \begin{center}
      \begin{tabular}{rr|rr|r}
        $k$& $n$& $T_{k}$ & $TH_{k}$ & \bigstrut\\
        \hline
        &&\multicolumn{2}{c|}{$\bdw\geq$} & $\bdw \leq$\\
        \hline
        \hline
        7 & 49 & 10 & 12 &  14 \\
        8 & 64 & 11 & 13 &  16 \\
        9 & 81 & 11 & 14 &  18 \\
        10 & 100& 14 & 15 & 20 \\
        \hline
      \end{tabular}
      \caption{Summary of bounds for the bandwidth.}
      \label{tab:sumBounds}
    \end{center}
  \end{table}

  \medskip
  We summarize the bandwidth information for all
  variations of the Torus graphs
  in Table \ref{tab:sumBounds}.
  Our partitioning approach provides nontrivial lower bounds on all instances.

  Now, let us provide some information on computation time.
  To compute 4-partitions for graphs with 49 vertices we need about $20$ seconds,
  for 5-partitions about $30$ seconds, and for 6-partitions about $90$ seconds.
  On the other hand, to compute a 4-partition (6-partition) on a graph
  with 100 vertices, our ADMM code needs about $200$ seconds ($700$ seconds).
  Clearly, computation times increase with respect to increasing  partition
  sizes and number of vertices of the graphs.
  However, we obtain bounds in reasonable time for all tested graphs.
  All experiments were performed on a Windows 7 64-bit machine equipped with
  an Intel Core i5-5300U (2$\times$2300 MHz) and    12 GB RAM using MATLAB 2016b.

  \paragraph{Hypercubes}
  Results for the  hypercubes $Q_5$, $Q_6$, and $Q_7$ are summarized in
  Table \ref{tab:Ham}.
  The table reads similar to the previous tables.
  To show a lower bound of 10 for $\bdw(Q_5)$, our ADMM needs only 4 seconds.
  For comparison purposes we computed a lower bound for $Q_5$ and the case $k=32$.
  Thus, we solved the QAP relaxation for that instance  and obtained $11$ as the
  lower bound of the BP.

  For the hypercube $Q_6$ the 4-partition with $m_2=m_3=17$ and $r=1$ yields a
  positive lower  bound, and therefore    $\bdw(Q_6)\geq 18$.
  We also compute the 6-partition with $m=\transpose{(15,9,8,9,8,15)}$ and $r=2$,
  and obtain a positive lower bound, which leads again to the conclusion that
   $\bdw(Q_6)\geq 18$.
  Finally,  we prove a lower bound of 33 for the hypercube $Q_7$.

  \begin{table}[!ht]
    \begin{center}
      \begin{tabular}{rrrrrr|rr|r|r}
        \multicolumn{10}{c}
        {Hypercube $Q_{5}  \qquad n=32, ~\bdw=13$ } \bigstrut\\
        \hline\hline
        $m_1$ & $m_2$ & $m_3$ &$m_4$ & - & - & $ub $ & $lb$ & $\bdw \geq$ & \\
        \hline\hline
        6 &10 &10 &  6 &   &    & 0 & - & & \\
        7 & 9 & 9 &  7 &   &    & 4 & 0.99 & 10 &    \\ \hline
        \multicolumn{10}{c}
        {Hypercube $Q_{6}  \qquad n=64, ~\bdw = 23$} \bigstrut\\
        \hline\hline
        $m_1$ & $m_2$ & $m_3$ &$m_4$ & $m_5$ & $m_6$ & $ub $ & $lb$ & $\bdw \geq$ &r \\
        \hline\hline
        15 &17 &17 & 15 &   &    & 10 & 1.18 & 18 & 1 \\
        15 & 9 & 8 &  9 & 8  &15 & 14 & 1.18 & 18 & 2   \\
        14 & 9 & 9 &  9 & 9  &14 & 9 & - & & \\
        \hline\hline
        \multicolumn{10}{c}
        {Hypercube $Q_{7}  \qquad n=128, ~\bdw = 43$} \bigstrut\\
        \hline\hline
        $m_1$ & $m_2$ & $m_3$ &$m_4$ & $m_5$ &  & $ub $ & $lb$ & $\bdw \geq$ & \\
        \hline\hline
        33 &31 &31 & 33 &   &    &19 & - & &    \\
        34 &30 &30 & 34 &   &    &31 & 3.11 & 31 &    \\ \hline
        16 &32 &32 & 32 &16 &    &18 & 0.93 & 33 &    \\
        \hline
      \end{tabular}
      \caption{Hypercubes.}
      \label{tab:Ham}
    \end{center}
  \end{table}

  \subsection{Bandwidth of Matrices from Applications}

  In this section, we evaluate the performance of our approach
  on  matrices that are given by real-world applications.
  We collected symmetric matrices, having $48$ to $115$ vertices.
  These are taken  from the HB, Pothen, and Pajek groups of the
  SuiteSparse Matrix Collection \cite{matrixAM}.
  We also selected matrices from the Newman collection
  available on the NIST Matrix Market \cite{matrixMarket}.

  Considering the \BandwidthProblem, only the structural properties
  of the matrices are of interest.
  Therefore,  for a matrix $A$, we set  $\diag(A)=0$.
  Moreover, we set all nonzero entries  equal to one.

  In our computational evaluation, we select the partitioning
  $m$ such that $m_2=\ldots=m_{k-1}$,
  $m_1= \lfloor \frac{n-d}{2} \rfloor$, and $m_k= \lceil \frac{n-d}{2} \rceil$
  where $d=\sum_{i=2}^{k-1} m_i$.
  We set $r=1$, except when applying the 6-partition to \textit{adjnoun} and \textit{football}
  where we had to set $r=2$. In the later case, we apply Corollary \ref{cor:alternatingM}.

  We summarize the results in Table \ref{tab:literature}.
  We provide the number of nodes (column labeled $n$) and the number of
  edges (column labeled $|E(G)|$). The column labeled $(bdw \leq)$ provides
  an upper bound on the bandwidth which we found by running a Simulated
  Annealing heuristic. We did not find any bandwidth information on these
  data in the literature. We also determined the density relative to the
  bandwidth, i.e., proportion of edges within the bandwidth, in the column labeled
  (bdw-dens). Finally, and most interestingly, we provide lower bounds based
  on $k$-partitions for $k \in \{3,4,5,6 \}$.
  The results in the column for $k=3$ reflect the previous state-of-the-art
  using 3-partitions. The remaining columns show the improvement of the lower
  bound using partitions into $k \in \{4,5,6\}$ blocks.
  The lower bound is substantially improved in all cases.
  These results clearly indicate that our general partition approach
  yields a significant improvement over the 3-partition bounds from
  \cite{JuMo:99,HelReMoPo:95,PoRe:07,Rendl2017}.

  \begin{table}[!ht]
    \begin{center}
      \begin{tabular}{lrrrrrrrr}
        \hline\hline
        Name  & $n$ & $|E(G)|$ & $\bdw~\leq$ & \bdw-dens &
        \multicolumn{4}{c}{$\bdw~\geq$}   \\
        &&&&& \multicolumn{4}{c}{partitioning}   \\
        &&&&& 3 & 4 & 5 & 6  \\

        \hline\hline
        DWT59   & 59  & 104 &  6 & 0.381 & 3 & 4 & 4 & 5\\
        DWT87   & 87  & 227 & 10 & 0.278 & 5& 6 & 7 & 8 \\
        \hline
        NOS4    & 100 & 247 & 10 & 0.261 & 6 & 7 & 7 & 8 \\
        \hline
        ASH85   &  85 & 219 &  9 & 0.304 & 4 & 6 & 7 & 7\\
        CAN61   &  61 & 248 & 13 & 0.353 & 5 & 9 & 9 & 11\\
        CAN73   &  73 & 152 & 16 & 0.147 & 7 & 11 & 14 & 14\\
        CAN96   &  96 & 336 & 13 & 0.290 & 7 & 10 & 11 & 12\\
        \hline
        GD97-b  &  47 & 132 & 15 & 0.226 & 5 & 11 & 12 & 11\\
        \hline
        mesh1e1 &  48 & 129 & 11 & 0.279 & 6 & 9 & 10 & 10 \\
        sphere2 &  66 & 192 & 13 & 0.250 & 7 & 9 & 11 & 12 \\
        \hline
        dolphins & 62 & 159 &  13 & 0.222 & 7 & 9 & 11 & 11 \\
        lesmis   & 77 & 254 &  20 & 0.191 & 5 & 11 & 16 & 17\\
        polbooks & 105 & 441 & 20 & 0.233 & 9 & 11 & 14  & 17  \\
        adjnoun  & 112 & 425 & 39 & 0.119 & 23 & 32 & 31 & 32\\
        football & 115 & 613 & 37 & 0.173 & 28 & 33 & 33 & 33 \\
        \hline
      \end{tabular}
      \caption{Graphs from the literature.}
      \label{tab:literature}
    \end{center}
  \end{table}

  \subsection{Discussion}

  Based on our computational experiments we reach the following conclusions.
  \begin{itemize}
    \item
    The partitioning approach leads to acceptable lower bounds for the \BandwidthProblem.
    Our results indicate that the bounds get weaker as the number of nodes
    increases. This should come as no surprise in view of the
    hardness results known for the \BandwidthProblem.

    \item
    Our approach offers some flexibility in choosing the number $k$ of partition
    blocks to estimate the bandwidth.
    A larger $k$ would result in tighter bounds at higher computational cost.

    \item
    Further tightening of the semidefinite models is possible by adding
    additional constraints, e.g., triangle inequalities.
    This results in SDPs which require a refined computational setup.

    \item
    We could prove significantly better lower bounds for the \BandwidthProblem compared to
    the previous state-of-the-art of using 3-partitions.
  \end{itemize}

  \section{Summary and Conclusion}
  We have shown that the partition approach provides a versatile tool to obtain
  lower bounds for the bandwidth of a graph.
  The choice of the model parameters $k$, $m$, and $r$ are highly
  problem dependent. However, our experiments indicate that
  even with a small number of partition blocks ($k \ll n$)
  we are able to derive nontrivial lower bounds on the bandwidth, even
  for very  sparse graphs.
  Further research is necessary to explore this approach for larger graphs.


\begin{thebibliography}{42}
    \providecommand{\natexlab}[1]{#1}
    \providecommand{\url}[1]{\texttt{#1}}
    \expandafter\ifx\csname urlstyle\endcsname\relax
    \providecommand{\doi}[1]{doi: #1}\else
    \providecommand{\doi}{doi: \begingroup \urlstyle{rm}\Url}\fi

    \bibitem[Balogh et~al.(2006)Balogh, Mubayi, and Pluh\'ar]{BALOGH200643}
    József Balogh, Dhruv Mubayi, and Andr\'as Pluh\'ar.
    \newblock On the edge-bandwidth of graph products.
    \newblock \emph{Theoretical Computer Science}, 359\penalty0 (1):\penalty0
    43--57, 2006.
    \newblock ISSN 0304-3975.
    \newblock \doi{10.1016/j.tcs.2006.01.046}.

    \bibitem[Berger-Wolf and Reingold(2002)]{BergWoRe:02}
    Tanya~Y. Berger-Wolf and Edward~M. Reingold.
    \newblock Index assignment for multichannel communication under failure.
    \newblock \emph{IEEE Transactions on Information Theory}, 48\penalty0
    (10):\penalty0 2656--2668, Oct 2002.
    \newblock ISSN 0018-9448.
    \newblock \doi{10.1109/TIT.2002.802643}.

    \bibitem[Blum et~al.(2000)Blum, Konjevod, Ravi, and Vempala]{BlKoRaVe}
    Avrim Blum, Goran Konjevod, R.~Ravi, and Santosh Vempala.
    \newblock Semi-definite relaxations for minimum bandwidth and other
    vertex-ordering problems.
    \newblock \emph{Theoretical Computer Science}, 235\penalty0 (1):\penalty0
    25--42, 2000.
    \newblock ISSN 0304-3975.
    \newblock \doi{10.1016/S0304-3975(99)00181-4}.

    \bibitem[Burkard and Rendl(1984)]{BURKARD1984169}
    Rainer~E. Burkard and Franz Rendl.
    \newblock A thermodynamically motivated simulation procedure for combinatorial
    optimization problems.
    \newblock \emph{European Journal of Operational Research}, 17\penalty0
    (2):\penalty0 169--174, 1984.
    \newblock ISSN 0377-2217.
    \newblock \doi{10.1016/0377-2217(84)90231-5}.

    \bibitem[Cavero et~al.(2021)Cavero, Pardo, Laguna, and
    Duarte]{CAVERO2021105116}
    Sergio Cavero, Eduardo~G. Pardo, Manuel Laguna, and Abraham Duarte.
    \newblock Multistart search for the cyclic cutwidth minimization problem.
    \newblock \emph{Computers \& Operations Research}, 126:\penalty0 105--116,
    2021.
    \newblock ISSN 0305-0548.
    \newblock \doi{doi.org/10.1016/j.cor.2020.105116}.

    \bibitem[Chinn et~al.()Chinn, Chv\'atalov\'a, Dewdney, and Gibbs]{ChChDeGi:82}
    P.~Z. Chinn, J.~Chv\'atalov\'a, A.~K. Dewdney, and N.~E. Gibbs.
    \newblock The bandwidth problem for graphs and matrices—a survey.
    \newblock \emph{Journal of Graph Theory}, 6\penalty0 (3):\penalty0 223--254.
    \newblock \doi{10.1002/jgt.3190060302}.

    \bibitem[Chv\'atal(1970)]{Chvatal1}
    V\'aclav Chv\'atal.
    \newblock A remark on a problem of {H}arary.
    \newblock \emph{Czechoslovak Mathematical Journal}, 20\penalty0 (1):\penalty0
    109--111, 1970.
    \newblock URL \url{http://eudml.org/doc/12520}.

    \bibitem[Chv\'atalov\'a(1975)]{Chvalatalova}
    Jarmila Chv\'atalov\'a.
    \newblock Optimal labelling of a product of two paths.
    \newblock \emph{Discrete Mathematics}, 11\penalty0 (3):\penalty0 249--253,
    1975.
    \newblock ISSN 0012-365X.
    \newblock \doi{10.1016/0012-365X(75)90039-4}.

    \bibitem[Cuthill and McKee(1969)]{CutMcKee}
    Elizabeth Cuthill and James McKee.
    \newblock Reducing the bandwidth of sparse symmetric matrices.
    \newblock In \emph{Proceedings of the 1969 24th National Conference}, ACM '69,
    pages 157--172, New York, NY, USA, 1969. ACM.
    \newblock \doi{10.1145/800195.805928}.

    \bibitem[de~Klerk and Sotirov(2012)]{dKSo:12}
    Etienne de~Klerk and Renata Sotirov.
    \newblock Improved semidefinite programming bounds for quadratic assignment
    problems with suitable symmetry.
    \newblock \emph{Mathematical Programming}, 133\penalty0 (1):\penalty0 75--91,
    Jun 2012.
    \newblock ISSN 1436-4646.
    \newblock \doi{10.1007/s10107-010-0411-5}.

    \bibitem[de~Klerk et~al.(2013)de~Klerk, E.-Nagy, and Sotirov]{dKlSoNa}
    Etienne de~Klerk, Marianna E.-Nagy, and Renata Sotirov.
    \newblock On semidefinite programming bounds for graph bandwidth.
    \newblock \emph{Optimization Methods and Software}, 28\penalty0 (3):\penalty0
    485--500, 2013.
    \newblock \doi{10.1080/10556788.2012.709856}.

    \bibitem[D\'{\i}az et~al.(2002)D\'{\i}az, Petit, and Serna]{DiPeSe:02}
    Josep D\'{\i}az, Jordi Petit, and Maria Serna.
    \newblock A survey of graph layout problems.
    \newblock \emph{ACM Comput. Surv.}, 34\penalty0 (3):\penalty0 313--356,
    September 2002.
    \newblock ISSN 0360-0300.
    \newblock \doi{10.1145/568522.568523}.

    \bibitem[Dunagan and Vempala(2001)]{DuVe:01}
    John Dunagan and Santosh Vempala.
    \newblock On euclidean embeddings and bandwidth minimization.
    \newblock In Michel Goemans, Klaus Jansen, Jos{\'e} D.~P. Rolim, and Luca
    Trevisan, editors, \emph{Approximation, Randomization, and Combinatorial
    Optimization: Algorithms and Techniques}, pages 229--240, Berlin, Heidelberg,
    2001. Springer Berlin Heidelberg.
    \newblock ISBN 978-3-540-44666-8.

    \bibitem[Garey et~al.(1978)Garey, Graham, Johnson, and Knuth]{GagrJoKn:78}
    Michael~R. Garey, R.~Graham, David Johnson, and D.~Knuth.
    \newblock Complexity results for bandwidth minimization.
    \newblock \emph{SIAM Journal on Applied Mathematics}, 34\penalty0 (3):\penalty0
    477--495, 1978.
    \newblock \doi{10.1137/0134037}.

    \bibitem[Haemers(1995)]{haemers}
    Willem~H. Haemers.
    \newblock Interlacing eigenvalues and graphs.
    \newblock \emph{Linear Algebra and its Applications}, 226--228:\penalty0
    593--616, 1995.
    \newblock ISSN 0024-3795.
    \newblock \doi{10.1016/0024-3795(95)00199-2}.
    \newblock Honoring J.J.Seidel.

    \bibitem[Harary(1967)]{Harary:67}
    Frank Harary.
    \newblock {Problem 16}.
    \newblock In Miroslav Fiedler, editor, \emph{Theory of graphs and its
    applications}, page 191. Czechoslovak Academy of Science, Prague, 1967.

    \bibitem[Harper(1964)]{Harper1964}
    Lawrence~H. Harper.
    \newblock {Optimal assignments of numbers to vertices}.
    \newblock \emph{Journal of the Society for Industrial and Applied Mathematics},
    12\penalty0 (1):\penalty0 131--135, 1964.

    \bibitem[Harper(1966)]{Harper66}
    Lawrence~H. Harper.
    \newblock Optimal numberings and isoperimetric problems on graphs.
    \newblock \emph{Journal of Combinatorial Theory}, 1\penalty0 (3):\penalty0
    385--393, 1966.
    \newblock ISSN 0021-9800.
    \newblock \doi{10.1016/S0021-9800(66)80059-5}.

    \bibitem[Helmberg et~al.(1995)Helmberg, Rendl, Mohar, and Poljak]{HelReMoPo:95}
    Christoph Helmberg, Franz Rendl, Bojan Mohar, and Svatopluk Poljak.
    \newblock A spectral approach to bandwidth and separator problems in graphs.
    \newblock \emph{Linear and Multilinear Algebra}, 39\penalty0 (1--2):\penalty0
    73--90, 1995.
    \newblock \doi{10.1080/03081089508818381}.

    \bibitem[Hochberg et~al.(1995)Hochberg, McDiarmid, and Saks]{HoMcDSa}
    Robert Hochberg, Colin McDiarmid, and Michael Saks.
    \newblock On the bandwidth of triangulated triangles.
    \newblock \emph{Discrete Mathematics}, 138\penalty0 (1):\penalty0 261--265,
    1995.
    \newblock ISSN 0012-365X.
    \newblock \doi{10.1016/0012-365X(94)00208-Z}.
    \newblock 14th British Combinatorial Conference.

    \bibitem[Hu and Sotirov(2019)]{HaoSotirov:19}
    Hao Hu and Renata Sotirov.
    \newblock On solving the quadratic shortest path problem.
    \newblock \emph{INFORMS Journal on Computing}, 32\penalty0 (2):\penalty0
    219--233, 2019.

    \bibitem[Hwang and Lagarias(1977)]{HwLa77}
    F.K. Hwang and J.C. Lagarias.
    \newblock Minimum range sequences of all k-subsets of a set.
    \newblock \emph{Discrete Mathematics}, 19\penalty0 (3):\penalty0 257--264,
    1977.
    \newblock ISSN 0012-365X.
    \newblock \doi{10.1016/0012-365X(77)90105-4}.

    \bibitem[Juvan and Mohar(1993)]{JuMo:99}
    Martin Juvan and Bojan Mohar.
    \newblock Laplace eigenvalues and bandwidth-type invariants of graphs.
    \newblock \emph{J. Graph Theory}, 17\penalty0 (3):\penalty0 393--407, July
    1993.
    \newblock ISSN 0364-9024.
    \newblock \doi{10.1002/jgt.3190170313}.

    \bibitem[Lai and Williams(1999)]{LaWi:99}
    Yung-Ling Lai and Kenneth Williams.
    \newblock A survey of solved problems and applications on bandwidth, edgesum,
    and profile of graphs.
    \newblock \emph{J. Graph Theory}, 31\penalty0 (2):\penalty0 75--94, June 1999.
    \newblock ISSN 0364-9024.

    \bibitem[Leung et~al.(1984)Leung, Vornberger, and
    Witthoff]{1fb351e298a049cc855e645215159c30}
    {Joseph Y.T.} Leung, Oliver Vornberger, and {James D.} Witthoff.
    \newblock On some variants of the bandwidth minimization problem.
    \newblock \emph{SIAM Journal on Computing}, 13\penalty0 (3):\penalty0 650--667,
    1984.
    \newblock ISSN 0097-5397.
    \newblock \doi{10.1137/0213040}.

    \bibitem[Li et~al.(1981)Li, Tao, and Shen]{LiTapShen}
    Q.~Li, M.Q. Tao, and Y.Q. Shen.
    \newblock The bandwidth of the discrete tori ${C}_m \times {C}_n$.
    \newblock \emph{J. China Univ. Sci. Tech.}, \penalty0 (11):\penalty0 1–16,
    1981.

    \bibitem[Mart\'{i} et~al.(2013)Mart\'{i}, Pantrigo, Duarte, and
    Pardo]{MARTI2013137}
    Rafael Mart\'{i}, Juan~J. Pantrigo, Abraham Duarte, and Eduardo~G. Pardo.
    \newblock Branch and bound for the cutwidth minimization problem.
    \newblock \emph{Computers \& Operations Research}, 40\penalty0 (1):\penalty0
    137--149, 2013.
    \newblock ISSN 0305-0548.
    \newblock \doi{doi.org/10.1016/j.cor.2012.05.016}.

    \bibitem[Monien(1986)]{Monien}
    Burkhard Monien.
    \newblock The bandwidth minimization problem for caterpillars with hair length
    3 is {NP}-complete.
    \newblock \emph{SIAM J. Algebraic Discrete Methods}, 7\penalty0 (4):\penalty0
    505--512, October 1986.
    \newblock ISSN 0196-5212.
    \newblock \doi{10.1137/0607057}.

    \bibitem[{NIST}()]{matrixMarket}
    {NIST}.
    \newblock {Matrix Market}.
    \newblock URL \url{http://math.nist.gov/MatrixMarket/}.
    \newblock [Online; accessed 29-Sept-2018].

    \bibitem[Oliveira et~al.(2018)Oliveira, Wolkowicz, and Xu]{HenryADMM}
    Danilo~Elias Oliveira, Henry Wolkowicz, and Yangyang Xu.
    \newblock {ADMM} for the {SDP} relaxation of the {QAP}.
    \newblock \emph{Mathematical Programming Computation}, 10\penalty0
    (4):\penalty0 631--658, Dec 2018.
    \newblock ISSN 1867-2957.
    \newblock \doi{10.1007/s12532-018-0148-3}.

    \bibitem[Papadimitriou(1976)]{Papadem:76}
    Christos~H. Papadimitriou.
    \newblock The {NP}-completeness of the bandwidth minimization problem.
    \newblock \emph{Computing}, 16\penalty0 (3):\penalty0 263--270, Sep 1976.
    \newblock ISSN 1436-5057.
    \newblock \doi{10.1007/BF02280884}.

    \bibitem[Povh and Rendl(2007)]{PoRe:07}
    Janez Povh and Franz Rendl.
    \newblock A copositive programming approach to graph partitioning.
    \newblock \emph{SIAM Journal on Optimization}, 18\penalty0 (1):\penalty0
    223--241, 2007.
    \newblock \doi{10.1137/050637467}.

    \bibitem[Rendl and Sotirov(2018)]{Rendl2017}
    Franz Rendl and Renata Sotirov.
    \newblock The min-cut and vertex separator problem.
    \newblock \emph{Computational Optimization and Applications}, 69\penalty0
    (1):\penalty0 159---187, Jan 2018.
    \newblock ISSN 1573-2894.
    \newblock \doi{10.1007/s10589-017-9943-4}.

    \bibitem[Rodriguez-Tello et~al.(2008)Rodriguez-Tello, Hao, and
    Torres-Jimenez]{RODRIGUEZTELLO20083331}
    Eduardo Rodriguez-Tello, Jin-Kao Hao, and Jose Torres-Jimenez.
    \newblock An effective two-stage simulated annealing algorithm for the minimum
    linear arrangement problem.
    \newblock \emph{Computers \& Operations Research}, 35\penalty0 (10):\penalty0
    3331--3346, 2008.
    \newblock ISSN 0305-0548.
    \newblock \doi{10.1016/j.cor.2007.03.001}.

    \bibitem[Rodriguez-Tello et~al.(2015)Rodriguez-Tello, Romero-Monsivais,
    Ramirez-Torres, and Lardeux]{RODRIGUEZTELLO201517}
    Eduardo Rodriguez-Tello, Hillel Romero-Monsivais, Gabriel Ramirez-Torres, and
    Fr\'{e}d\'{e}ric Lardeux.
    \newblock Tabu search for the cyclic bandwidth problem.
    \newblock \emph{Computers \& Operations Research}, 57:\penalty0 17--32, 2015.
    \newblock \doi{10.1016/j.cor.2014.11.013}.

    \bibitem[Smithline(1995)]{Smith}
    Lawren Smithline.
    \newblock Bandwidth of the complete k-ary tree.
    \newblock \emph{Discrete Mathematics}, 142\penalty0 (1):\penalty0 203--212,
    1995.
    \newblock ISSN 0012-365X.
    \newblock \doi{10.1016/0012-365X(93)E0219-T}.

    \bibitem[{Tim Davis}()]{matrixAM}
    {Tim Davis}.
    \newblock {The SuiteSparse Matrix Collection}.
    \newblock URL \url{https://sparse.tamu.edu}.
    \newblock [Online; accessed 11-Nov-2019].

    \bibitem[Turner(1986)]{Turner86}
    Jonathan~S. Turner.
    \newblock On the probable performance of heuristics for bandwidth minimization.
    \newblock \emph{SIAM Journal on Computing}, 15\penalty0 (2):\penalty0 561--580,
    1986.
    \newblock \doi{10.1137/0215041}.

    \bibitem[Unger(1998)]{Unger1998}
    Walter Unger.
    \newblock {The complexity of the approximation of the bandwidth problem}.
    \newblock In \emph{Proceedings 39th Annual Symposium on Foundations of Computer
    Science (Cat. No.98CB36280)}, pages 82--91, Nov 1998.
    \newblock \doi{10.1109/SFCS.1998.743431}.

    \bibitem[van Dam and Sotirov(2015)]{DamSotirov15}
    Edwin~R. van Dam and Renata Sotirov.
    \newblock On bounding the bandwidth of graphs with symmetry.
    \newblock \emph{INFORMS Journal on Computing}, 27\penalty0 (1):\penalty0
    75--88, 2015.
    \newblock \doi{10.1287/ijoc.2014.0611}.

    \bibitem[Wolkowicz and Zhao(1999)]{WOLKOWICZ1999461}
    Henry Wolkowicz and Qing Zhao.
    \newblock {Semidefinite programming relaxations for the graph partitioning
    problem}.
    \newblock \emph{Discrete Applied Mathematics}, 96\penalty0 (Supplement
    C):\penalty0 461--479, 1999.
    \newblock ISSN 0166-218X.
    \newblock \doi{10.1016/S0166-218X(99)00102-X}.

    \bibitem[Zhao et~al.(1998)Zhao, Karisch, Rendl, and Wolkowicz]{ZhKaReWo:98}
    Qing Zhao, Stefan~E. Karisch, Franz Rendl, and Henry Wolkowicz.
    \newblock Semidefinite programming relaxations for the quadratic assignment
    problem.
    \newblock \emph{Journal of Combinatorial Optimization}, 2\penalty0
    (1):\penalty0 71--109, Mar 1998.
    \newblock ISSN 1573-2886.
    \newblock \doi{10.1023/A:1009795911987}.

  \end{thebibliography}
  \end{document}